\documentclass[11pt,russian]{article}   
\usepackage[russian]{babel}             
\normalfont

\usepackage{amssymb,amsmath}
\usepackage[dvips]{color,graphicx}

\newtheorem{theorem}{Теорема}
\newtheorem{lemma}{Лемма}

\begin{document}

\noindent
\textbf{\LARGE Восстановление собственных функций \\ $q$-ичного $n$-мерного гиперкуба
\footnote{Исследование выполнено за счет гранта Российского научного фонда (проект \No 14-11-00555)}}
\date{}

\vspace*{5mm} \noindent
\textsc{А.Ю. Васильева} \hfill \texttt{vasilan@math.nsc.ru} \\
{\small Институт математики им. С.Л. Соболева СО РАН,  \\
Новосибирский Государственный Университет, Новосибирск, Россия} \\

\medskip

\begin{center}
\parbox{11,8cm}{\footnotesize
Установлено, что значения произвольной собственной функции $q$-ичного $n$-мерного гиперкуба могут быть однозначно восстановлены во всех вершинах шара, если известны ее значения на соответствующей сфере; в терминах собственного числа и радиуса шара сформулированы достаточные условия для такого восстановления. Показано, что в случае, когда значения собственной функции заданы на сфере радиуса, равного номеру собственного числа, однозначно восстанавливаются все значения собственной функции; аналогично предыдущему случаю, выписаны достаточные числовые условия. }
\end{center}

\baselineskip=0.9\normalbaselineskip


Предметом исследования являются собственные функции $n$-мерного $q$-ичного гиперкуба, 
задача состоит в однозначном, полном или частичном, восстановлении функции по части ее значений, а главным инструментом является формула связи локальных распределений собственной функции в двух ортогональных гранях \cite{VasOC}. Задача восстановления ранее рассматривалась в двоичном случае, например в \cite{OHed,Vas2013}. Еще раньше, в \cite{Avg95}, было доказано, 
что 1-совершенный код длины $n$
однозначно определяется своими вершинами на среднем слое гиперкуба, 
что и положило начало для изучения подобных вопросов 
для более широких классов объектов в более широких классах графов.

Статья организована следующим образом. Раздел \ref{Prelim} является вспомогательным. В разделе \ref{LocDistr} приводится формула связи локальных распределений собственной функции в паре ортогональных граней и дается явное выражение компонент распределения в одной грани через компоненты в другой грани в случае различных соотношений размерностей этих граней с собственным числом функции (леммы \ref{vjort}, \ref{opjat}). Эти леммы применяются в разделе \ref{PartRec} для доказательства возможности, при некоторых числовых условиях, восстановления собственной функции в шаре по ее значениям на соответствующей сфере, формулируется основная теорема \ref{main}. В разделе \ref{FullRec} устанавливается, при аналогичных числовых условиях, что собственная функция из $h$-го собственного подпространства, $h=0,\ldots,n,$ полностью определяется своими значениями на сфере радиуса $h$.

\section{Предварительные сведения}\label{Prelim}

Пусть $q>2$ -- целое положительное число, не обязательно являющееся степенью простого.
Множество ${\bf F}_q=\{0,1,\ldots,q-1\}$ будем рассматривать как группу по сложению по модулю $q$, а
${\bf F}_q^n$ -- как абелеву группу
${\bf F}_q\times\ldots\times {\bf F}_q$.
Мы изучаем функции, заданные на множестве вершин графа ${\bf F}_q^n$ $q$-ичного $n$-мерного гиперкуба, множество ребер которого определяется как множество пар $q$-ичных наборов, различающихся в точности в одной позиции.

Везде далее через $I$ будем обозначать подмножество множества $\{1,\ldots, n\}$, а через $\overline{I}$ -- его дополнение:  $\overline{I}=\{1,\ldots,n\}\backslash I$.
Пусть $\alpha\in{\bf F}_q^n$ -- произвольная вершина гиперкуба. Обозначим \emph{носитель} вершины $\alpha$, т.е. множество ее ненулевых позиций, через $s(\alpha)$; мощность носителя -- это \emph{вес Хэмминга} вершины $\alpha$, он обозначается через $wt(\alpha)$;
\emph{расстояние Хэмминга} $\rho(\alpha,\beta)$ между двумя вершинами $\alpha$ и $\beta$  -- это вес Хэмминга их разности $\alpha-\beta$;
будем обозначать через $W_i(\alpha)$ множество всех вершин на расстоянии $i$ от $\alpha$, т.е. \emph{сферу} радиуса $i$ с центром $\alpha$, а через $B_i(\alpha)$ -- множество вершин на расстоянии не более $i$ от $\alpha$, т.е. \emph{шар} радиуса $i$ с центром $\alpha$.
Будем обозначать через $W_i$ и $\Gamma^I$ множества $W_i(\alpha)$ и $\Gamma^I(\alpha)$ в том случае, когда $\alpha= (0,\ldots,0)$.
Положим
$$\Gamma^I(\alpha)=\{\beta\in{\bf F}_q^n \ : \ \beta_i=\alpha_i \ \forall \ i\notin I\} ,$$
тогда $\Gamma^I(\alpha)$ называется $|I|$-\emph{мерной гранью} (или \emph{подкубом}),
как нетрудно заметить, она имеет структуру гиперкуба меньшей размерности ${\bf F}_q^{|I|}$.
Две грани $\Gamma^I(\alpha)$ и $\Gamma^J(\beta)$ называются \emph{ортогональными}, если $\overline{J}=I$. Из определения сразу следует, что две ортогональные грани имеют в точности одну общую вершину.
Для произвольных $\alpha,\beta\in {\bf F}_q^n$
обозначим $\langle\alpha,\beta\rangle=\alpha_1\beta_1+\ldots+\alpha_n\beta_n(\mod q) .$

Рассмотрим пространство всех комплекснозначных функций, заданных на вершинах $q$-ичного $n$-мерного гиперкуба:
$$V= \left\{ f: {\bf F}_q^n\longrightarrow \mathbb{C}\right\} .$$
Функцию можно рассматривать как $q^n$-мерный вектор ее значений.

\emph{Схема отношений Хэмминга} (введение в теорию схем отношений можно найти, например, в \cite{Dels})
состоит из основного множества ${\bf F}_q^n$ с заданными на нем $n+1$ отношениями $R_0,\ldots,R_n$, определяемыми по расстоянию Хэмминга:
для произвольных $\alpha,\beta\in {\bf F}_q^n$ и $i=0,1,\ldots,n$ верно ;$$(\alpha,\beta)\in R_i
 \ \Leftrightarrow \   \rho(\alpha,\beta)=i.$$

Пусть $D_i=D^{q,n}_i$ -- матрица смежности $i$-го отношения $R_i$, т.е. квадратная матрица
порядка $q^n$ с элементами
$$(D^{q,n}_i)_{\alpha,\beta} =
\left\{\begin{array}{ll} 1, & \rho(\alpha,\beta)=i , \\ 0, & \rho(\alpha,\beta)\neq i . \end{array}\right.$$
Отметим, что $D=D_1=D^{q,n}_1$ является матрицей смежности графа гиперкуба ${\bf F}_q^n$.
Комплексное число $\lambda$ будем называть собственным значением графа, если оно является собственным значением его матрицы смежности.
В частности, $\lambda$ -- ,
если оно является собственным числом матрицы $D$.
Известно \cite{Dels}, что собственные значения графа $q$-значного $n$-мерного гиперкуба равны
$$(q-1)n-qh,  \ \ \ h=0,1,\ldots,n ,$$
в этом случае число $h$ называется \emph{номером собственного значения} $\lambda_h$.
Другими словами, номер собственного значения $\lambda$ равен
$h=h(\lambda)= \frac{(q-1)n-\lambda}{q}$.
Соответствующие собственные функции (мы будем называть их \emph{$\lambda$-функциями}) удовлетворяют уравнениям
\begin{equation} \label{ball}  \sum_{\beta\in W_1(\alpha)} f(\beta) = \lambda f(\alpha), \ \ \ \alpha\in{\bf F}_q^n,
\end{equation}
или, обозначая через $f$ вектор значений функции $f$:
$$Df=\lambda f .$$
Все матрицы смежности $D^{q,n}_i, \ i=0,\ldots,n,$ \ $q$-ичной $n$-мерной схемы Хэмминга
имеют общую систему ортогональных собственных подпространств $V_h, \ h=0,\ldots,n,$
а собственное значение матрицы $D^{q,n}_i$ на $V_h$ есть $P_i^{(q)}(h;n)$, где
$$P^{(q)}_i(t;N)= \sum_{j=0}^i (-1)^j (q-1)^{i-j}\binom{t}{j}\binom{N-t}{i-j} - $$
многочлен Кравчука от переменной $t$. Нужные нам значения многочленов Кравчука
могут быть заданы как коэффициенты многочлена от переменных $x$ и $y$:
$$(x-y)^{t} (x+(q-1)y)^{N-t} = \sum_{i=0}^N P^{(q)}_i(t;N)y^i x^{N-i} . $$
Второй базис пространства матриц
$\{\sum_{i=0}^n a_i D_i, \ \ a_i\in \mathbb{C}, \ i=0,\ldots,n\}$
дают так называемые примитивные идемпотенты $J_0,\ldots,J_n$, задающие проектирование на подпространства $V_0,\ldots,V_n$ соответственно.
Связь двух базисов описывается опять же при помощи многочленов Кравчука:
\begin{equation}\label{DinJ}
D_i=\sum_{j=0}^k P_i^{(q)}(j;n) J_j .
\end{equation}
Полезным свойством базиса из примитивных идемпотентов является то, что если
$A=\sum_{i=0}^n a_i J_i$, то $A^{-1}=\sum_{i=0}^n a_i^{-1} J_i$.

\section{Локальные распределения}\label{LocDistr}
Как обобщение хорошо известного понятия \emph{весового распределения} функции
$v^{f}(\alpha)=(v^{f}_0(\alpha),\ldots,v^{f}_{n}(\alpha))$, где
$$v^{f}_j (\alpha)= \sum_{\beta\in W_j(\alpha)} f(\beta) , \ \ \ i=0,\ldots,n,$$
для произвольной функции $f$ относительно произвольной вершины $\alpha\in {\bf F}_q^n$,
введем определение локального распределения.
Положим
$$v^{I,f}_j (\alpha)= \sum_{\beta\in\Gamma^I(\alpha)\bigcap W_j(\alpha)} f(\beta) ,
 \ \ \ i=0,\ldots,|I|,$$
вектор $v^{I,f}(\alpha)=(v^{I,f}_0(\alpha),\ldots,v^{I,f}_{|I|}(\alpha))$ называется \emph{локальным распределением} функции $f$ в грани $\Gamma^I(\alpha)$ относительно вершины $\alpha$, или короче, $(I,\alpha)$-\emph{локальным распределением} $f$.  Многочлен $g^{I,\alpha}_f(x,y)$ называется $(I,\alpha)$-\emph{ло\-каль\-ным энумератором} функции $f$  в грани $\Gamma^I(\alpha)$ относительно вершины $\alpha$, или $(I,\alpha)$-\emph{локальным энумератором} $f$, если
\begin{equation}\label{loc-enum}
g^{I,\alpha}_f(x,y) = \sum_{j=0}^{|I|} v_j^{I,f}(\alpha) y^j x^{|I|-j} =
\sum_{\beta\in\Gamma^I(\alpha)} f(\beta) y^{|s(\beta)|}
x^{|I|-|s(\beta)|} .
\end{equation}
Индекс $\alpha$ будем опускать во всех обозначениях, если $\alpha=(0,\ldots,0)$.

Взаимосвязь локальных энумераторов произвольной $\lambda$-функции в двух ортогональных гранях
дается следующей теоремой:
\begin{theorem}\cite{VasOC} \label{q-eig_loc}
Пусть $\lambda$ -- собственное число ${\bf F}_q^n$ с номером $h$, \ а $f$ -- произвольная $\lambda$-функция. Тогда для любого $\alpha\in{\bf F}_q^n$ выполнено
\begin{equation}\label{fla-loc}
(x+(q-1)y)^{h-|\overline{I}|} g^{\overline{I},\alpha}_f(x,y) = (x'+(q-1)y')^{h-|I|} g^{I,\alpha}_f(x',y') ,
\end{equation}
где $x'=x+(q-2)y$, и $y'=-y $.
\end{theorem}

Цель этого раздела -- при помощи Теоремы \ref{q-eig_loc} выписать явные формулы, выражающие компоненты $(\overline{I},\alpha)$-локального распределения произвольной $\lambda$-функции через компоненты ее $(I,\alpha)$-локального распределения. Как нетрудно увидеть в (\ref{fla-loc}), такие формулы будут иметь различный вид в зависимости от соотношения размерности грани $k=|I|$ и номера $h$ собственного значения гиперкуба.
Возможны следующие варианты:

I) $k\leq min\{h,n-h\}$,

II) $h<k\leq n-h$,

III) $n-h<k\leq h$,

IV) $k> max\{h,n-h\}$.

Далее подробно остановимся на случаях I) и III). Такой выбор станет очевиден в последнем разделе, сейчас скажем лишь, что явные формулы нужны нам для пошагового восстановления $\lambda$-функций в сферах, причем максимальный требующийся радиус равен $h$.

Сначала представим $g_f^{I,\alpha}(x',y')$
как многочлен от переменных $x, y$:
\begin{lemma}\label{l-g-polin}
Для произвольной функции $f: {\bf F}_q^n \rightarrow \mathbb{C}$ выполнено
\begin{equation}\label{g-polin}
g_f^{I,\alpha}(x+(q-2)y,-y)=\sum_{l=0}^k y^l x^{k-l} \sum_{i=0}^l (-1)^i v_i^{I,f}(\alpha) (q-2)^{l-i}\left(\begin{array}{c} k-i \\ l-i \end{array}\right) .
\end{equation}
\end{lemma}

{\bf Доказательство}.
По определению локального энумератора имеем:
$$g_f^{I,\alpha}(x+(q-2)y,-y)=\sum_{i=0}^k v_i^{I,f}(\alpha) (-y)^i (x+(q-2)y)^{k-i}$$
$$ =\sum_{i=0}^k v_i^{I,f}(\alpha) (-y)^i \sum_{l=0}^{k-i} (q-2)^l y^l x^{k-i-l} .$$
Меняя порядок суммирования, получим (\ref{g-polin}).
\hfill $\square$

Рассмотрим случай I).
Поскольку $k<h$ и $k<n-h$, то в формуле
\begin{equation} \label{loc-enum}
g_f^{\overline{I},\alpha}(x,y) =
(x-y)^{h-|I|} (x+(q-1)y)^{|\overline{I}|-h} g_f^{I,\alpha}(x+(q-2)y,-y) ,
\end{equation}
полученной из (\ref{fla-loc}), в правой части стоит произведение многочленов,
причем $(x-y)^{h-|I|} (x+(q-1)y)^{|\overline{I}|-h}$ есть производящая функция
значений многочленов Кравчука, а $g_f^{I,\alpha}(x+(q-2)y,-y)$ дается Леммой \ref{g-polin}.
Перемножая их, получаем:
\begin{lemma} \label{vjort}
Пусть $\lambda$ -- собственное число ${\bf F}_q^n$ с номером $h$, \ а $f$ -- произвольная $\lambda$-функция. Пусть $|I|=k\leq  min\{h,n-h\}$. Тогда для любого $j=0,\ldots,n-k$ выполнено
\begin{equation}
 v_j^{\overline{I}}(\alpha)= \sum_{i=0}^j r_{ij} v_i^I(\alpha), \ \ \ \  \mbox{где}
\end{equation}
\begin{equation}\label{rij}
r_{ij}=(-1)^i \sum_{l=0}^{j-i} P_{j-i-l}^{(q)}(h-k;n-2k)(q-2)^l \left(\begin{array}{c} k-i \\ l \end{array}\right)
\end{equation}
\end{lemma}


Перейдем к рассмотрению случая III),
пусть теперь $n-h<k\leq h$.
Перепишем формулу (\ref{loc-enum}):
\begin{equation} \label{n-h<h}
(x+(q-1)y)^{h+k-n} g^{\overline{I},\alpha}(x,y) = (x-y)^{h-k} g^{I,\alpha}(x+(q-2)y,-y)
\end{equation}
Учитывая, что здесь опять $(x+(q-1)y)^{h+k-n}$ и $(x-y)^{h-k}$ являются многочленами,
вычислим обе части равенства.

Во-первых,
$$(x+(q-1)y)^{h+k-n} g^{\overline{I},\alpha}(x,y) =
\sum_{j=0}^h y^j x^{h-j} \sum_{i=0}^j v_i^{\overline{I},f}(\alpha) (q-1)^{j-i}
\left(\begin{array}{c} h+k-n \\ j-i \end{array}\right) .$$
Во-вторых, применяя Лемму \ref{l-g-polin}, имеем:
$$(x-y)^{h-k} g_f^{I,\alpha}(x+(q-2)y,-y) = $$
$$\sum_{j=0}^h y^j x^{h-j} \sum_{l=0}^j (-1)^{j-l}\binom{h-k}{j-l}
 \sum_{i=0}^l (-1)^i v_i^{I,f}(\alpha) (q-2)^{l-i} \binom{k-i}{l-i} =$$
Меняя порядок суммирования и сравнивая внутреннюю сумму с определением многочленов Кравчука, получаем:
$$\sum_{j=0}^h y^j x^{h-j} \sum_{i=0}^j (-1)^i v_i^{I,f}(\alpha)
\sum_{l=i}^j (-1)^{j-l}\binom{h-k}{j-l} (q-2)^{l-i} \binom{k-i}{l-i} =$$
$$\sum_{j=0}^h y^j x^{h-j} \sum_{i=0}^j (-1)^i v_i^{I,f}(\alpha) P_{j-i}^{(q-1)}(h-k,h-i) .$$
Для каждого $j=0,\ldots,n-k<h$ приравниваем коэффициенты при $y^j x^{h-j}$
в левой и правой частях~(\ref{n-h<h}):
$$\sum_{i=0}^j v_i^{\overline{I},f}(\alpha) (q-1)^{j-i}
\left(\begin{array}{c} h+k-n \\ j-i \end{array}\right) =
\sum_{i=0}^j (-1)^i v_i^{I,f}(\alpha) P_{j-i}^{(q-1)}(h-k,h-i) .$$
Это система из $n-k+1$ уравнения относительно $n-k+1$ неизвестных
$v_0^{\overline{I},f}(\alpha),\ldots,v_{n-k+1}^{\overline{I},f}(\alpha)$
с квадратной матрицей $U$ с элементами
\begin{equation}\label{uji}
u_{ji} = (q-1)^{j-i} \left(\begin{array}{c} h+k-n \\ s-i \end{array}\right)
\ \ \ \ i,j=0,\ldots,n-k.
\end{equation}
Эта матрица -- нижнетреугольная с ненулевыми элементами на диагонали, а потому существует обратная к ней матрица $U^{-1} = (u^{\prime}_{si})_{s,i=0\ldots,j}$, тоже являющаяся нижнетреугольной.
Наша система разрешима и ее единственное решение имеет следующий вид, $j=0,\ldots,n-k$:
$$v_j^{\overline{I},f}(\alpha) =
\sum_{s=0}^j  u'_{js} \sum_{i=0}^s (-1)^i v_i^{I,f}(\alpha) P_{s-i}^{(q-1)}(h-k;h-i) .$$
Опять меняя порядок суммирования, получаем:
\begin{lemma} \label{opjat}
Пусть $\lambda$ -- собственное число ${\bf F}_q^n$ с номером $h$, \ а $f$ -- произвольная $\lambda$-функция. Пусть $k=|I|$ и $h-h<k\leq h$. Тогда для любого $j=0,\ldots,h-k$ выполнено
\begin{equation}\label{vi<min}
v_j^{\overline{I},f}(\alpha) = \sum_{i=0}^j  r_{ij}^k v_i^{I,f}(\alpha) , \ \ \ \ \mbox{где}
\end{equation}
\begin{equation}\label{rijk>min}
r_{ij}^k = (-1)^i \sum_{s=i}^j u'_{js} P_{s-i}^{(q-1)}(h-k;h-i) .
\end{equation}
а $u^{\prime}_{js}$ -- элементы матрицы $U^{-1}$, обратной к определенной в (\ref{uji}).
\end{lemma}
Итак, Леммы \ref{vjort} и \ref{opjat} дают нам линейное выражение компонент $(\overline{I},\alpha)$-ло\-каль\-но\-го распределения через компоненты $(I,\alpha)$-локального распределения, причем эти выражения зависят от  коэффициентов $r_{ij}^k$, вид которых зависит от $k$:
они определяются формулой (\ref{rij}) в случае $k\leq min\{n-h,h\}$
и формулой (\ref{rijk>min}) в случае $n-h<k\leq h$.

В заключение этого раздела отметим, что в отличие от случая, когда распределения рассматриваются относительно нулевой вершины, в случае произвольной другой вершины $\alpha$ гиперкуба компоненты локальных распределений относительно нее вычисляются суммированием по вершинам, вообще говоря, разного веса. Точнее, пусть $\alpha\in F_q^n$ имеет вес $k$. Положим $I=s(\alpha)$ и рассмотрим грань $\Gamma^I(\alpha)$. Это $k$-мерная грань, содержащая нулевую вершину.
Тогда
\begin{equation}\label{1v-sigma-delta}
v_i^{I,f}(\alpha)= \sigma_i^{I,f}(\alpha) + \delta_i^{I,f}(\alpha)  , \ \ \ \ \mbox{где}
\end{equation}
$$\sigma_i^{I,f}(\alpha) = \sum_{\beta\in W_k\bigcap\Gamma^I(\alpha)\bigcap W_i(\alpha)} f(\beta)
 \ \ \ \mbox{и} \ \ \
 \delta_i^{I,f}(\alpha) =
 \sum_{\beta\in B_{k-1}\bigcap\Gamma^I(\alpha)\bigcap W_i(\alpha)} f(\beta) ,$$
т.е. компонента $v_i^{I,f}(\alpha), \ i=0,\ldots,|I|,$
может быть представлена как сумма  $\sigma_i^{I,f}(\alpha)$ значений функции
по всем вершинам веса $k$
плюс сумма $\delta_i^{I,f}(\alpha)$ значений функции по вершинам меньшего веса.

\section{Частичное восстановление $\lambda$-функций}\label{PartRec}
Пусть $\lambda$ -- собственное значение гиперкуба ${\bf F}_q^n$ с номером $h=((q-1)n-\lambda)/q$,
а $f$ соответствующая ему $\lambda$-функция, и пусть $d\leq h$.

Зададимся следующим вопросом.
Предположим, что нам известны значения $f(\alpha)$ для всех вершин $\alpha$ веса Хэмминга $d$.
Возможно ли однозначно определить значения $f(\alpha)$ для всех вершин $\alpha$ с весом Хэмминга, меньшим, чем $d$?

Во-первых, вследствие дистанционной регулярности нашего графа весовое распределение произвольной собственной функции однозначно определяется номером собственного значения и значением функции в начальной вершине, а потому
$$\sum_{\alpha\in W_d} f(\alpha) = P^{(q)}_d(h;n) f({\bf 0}) ,$$
Эта формула следует из того, что в схеме отношений Хэмминга (см., напр.,\cite{Dels}):
$$ D_d f = P^{(q)}_d(h;n) f .$$
Применительно к нашей ситуации это означает,
что мы можем вычислить значение функции в нулевой вершине при условии, что
$P^{(q)}_d(h;n)\neq 0$ :
\begin{equation}\label{fo}
f({\bf 0}) = \frac{\sum_{\alpha\in W_d} f(\alpha)}{P^{(q)}_d(h;n)} .
\end{equation}
Далее попытаемся вычислить значения нашей функции в вершинах веса 1, затем 2, 3 и т.д,
используя индукцию по весу вершины (базис индукции мы только что получили).

Введем дополнительные обозначения: $S^I$ -- множество всех вершин гиперкуба,
имеющих носитель $I$;  \ $F^I$ -- вектор всех значений $f(\alpha), \ \alpha\in U^I$; \
$\Phi^I$ -- вектор всех компонент $v_{d-k}^{\overline{I},f}(\alpha), \ \alpha\in U^I$.
\begin{lemma} \label{insupp}
Пусть $\lambda$ -- собственное число ${\bf F}_q^n$ с номером $h$, \ а $f$ -- произвольная $\lambda$-функция. Пусть $|I|=k\leq h$. Тогда
\begin{equation}\label{inUI}
\sum_{i=0}^k r_{i,d-k}^k D_i^{q-1,k} F^I = \Phi^I - \Psi^I,
\end{equation}
где $D_i^{q-1,k}, \ i=0,\ldots,k,$ -- матрицы смежности $(q-1)$-ичной $k$-мерной схемы Хэмминга,
а вектор $\Psi^I$ зависит только от значений функции $f$ в вершинах веса меньше $k$
и не зависит от $F^I$.
\end{lemma}

{\bf Доказательство}.
Рассмотрим произвольную вершину $\alpha\in U^I$. При помощи лемм \ref{vjort} и \ref{opjat}
и свойства (\ref{1v-sigma-delta}) получаем:
$$v_{d-k}^{\overline{I}}(\alpha) =
\sum_{i=0}^{d-k} r_{i,d-k}^k v_i^{I,f}(\alpha) =
\sum_{i=0}^{d-k} r_{i,d-k}^k \left(\sigma_i^{I,f}(\alpha) + \delta_i^{I,f}(\alpha)\right) .$$
Это означает, что
\begin{equation}\label{Phi-Psi}
\sum_{i=0}^{d-k} r_{i,d-k}^k
\sum_{\beta\in S^I\bigcap W_i(\alpha)} f(\beta)
= v_{d-k}^{\overline{I}}(\alpha) -
    \sum_{i=0}^{d-k} r_{i,d-k} \sigma_i^-(\alpha) .
\end{equation}
Множество $U^I$ с отношениями по расстоянию на нем имеет структуру $(q-1)$-ичной $k$-мерной схемы Хэмминга, а потому совокупность соотношений (\ref{Phi-Psi}) для всех $\alpha\in S^I$
дает нам (\ref{inUI}).
\hfill $\square$

Перейдем к первой из двух основных теорем. Она позволит нам восстанавливать $\lambda$-функцию в шаре, зная ее значения на соответствующей сфере при некоторых дополнительных числовых условиях.

\begin{theorem} \label{main}
Пусть $\lambda$ -- собственное значение ${\bf F}_q^n$ с номером $h$, \
а $f$ -- произвольная $\lambda$-функция.
Пусть $d\leq h$  и
$\varphi: W_d\longrightarrow \mathbb{C}$ -- произвольная функция,
причем для всех $\alpha\in W_d$ выполнено $f(\alpha)=\varphi(\alpha) .$
Тогда для произвольного $\alpha\in B_d$ 
значение $f(\alpha)$ однозначно определено по функции $\varphi$
 при условии, что
\begin{equation}\label{maincond}
\sum_{i=0}^k r^k_{i, d-k} P^{(q-1)}_i(l,k) \neq 0 , \ \ \ \ k=0,\ldots,h, \ \ l=0,\ldots,k,
\end{equation}
где коэффициенты $r^k_{ij}$ определяются формулой (\ref{rij}) в случае $k\leq n-h$
и формулой (\ref{rijk>min}) в случае $k>n-h$.
\end{theorem}

{\bf Доказательство} проведем индукцией по весу вершин,
ее базис описан в начале этого раздела.

Теперь предположим, что (\ref{maincond}) выполнено
и значения функции $f$ в вершинах веса не более $k-1$
определяются однозначно по функции $\varphi$.
Пусть $I$ -- произвольное подмножество мощности  $k$ из $\{ 1,\ldots,n\}$
и вершина $\alpha$ имеет носитель $s(\alpha)=I$.
Очевидно, что носители всех вершин из грани $\Gamma^I(\alpha)$ содержатся в $I$.
По лемме \ref{insupp} имеем
\begin{equation}\label{M-FI}
M^{d,k} F^I = \Phi^I - \Psi^I, \ \ \ \ \mbox{где}
\end{equation}
\begin{equation}\label{matr}
M^{d,k}= \sum_{i=0}^k r^k_{i,d-k} D_i^{q-1,k}
\end{equation}
Согласно (\ref{DinJ}), представим матрицы смежности схемы Хэмминга
как линейные комбинации ее примитивных идемпотентов $J_l^{q-1,k}, \ l=0,\ldots,k$.
Тогда
$$M^{d,k}=  \sum_{i=0}^k r^k_{i,d-k} \sum_{l=0}^k P_i^{(q-1)}(l,k) J_l^{q-1,k} =
\sum_{l=0}^kJ_l^{q-1,k} \sum_{i=0}^k r^k_{i, d-k} P^{(q-1)}_i(l,k).$$
По условию (\ref{maincond}), все коэффициенты этой линейной комбинации не равны нулю, а потому матрица $M^{d,k}$ невырождена (это следует из свойств схем отношений) и система (\ref{M-FI})
имеет единственное решение. Перебирая различные подмножества $I$ мощности $k$, будем получать системы вида  (\ref{M-FI}), отличающиеся лишь правой частью, что не влияет на их однозначную разрешимость.
Таким образом, значения $f$ во всех вершинах веса $k$определяются однозначно по функции $\varphi$.
\hfill $\square$

\section{Полное восстановление $\lambda$-функций}\label{FullRec}
В этом разделе будет доказано, что все значения $\lambda$-функции однозначно определяются по ее
значениям на множестве вершин веса $h$, где $h$ -- номер собственного значения $\lambda$. Как и в предыдущей теореме, здесь будут фигурировать дополнительные числовые условия.

Как мы выяснили ранее, по значениям $\lambda$-функции на множестве вершин веса $h$ можно восстановить ее значения в вершинах меньшего веса. Таким образом, сейчас нам остается показать, что по значениям во всех вершинах шара радиуса $h$ однозначно восстанавливается вся функция.

\begin{lemma} \label{ball-cub}
Пусть $\lambda$ -- собственное значение ${\bf F}_q^n$ с номером $h$ и $f$ -- $\lambda$-функция. Пусть $\varphi: B_h\longrightarrow C$ -- произвольная функция. Предположим, что для всех $\alpha\in B_h$ выполнено $f(\alpha)=\varphi(\alpha) .$
Тогда все значения функции $f$ определены однозначно функцией $\varphi$.
\end{lemma}

Любая функция, заданная на гиперкубе, однозначно определяется своими коэффициентами Фурье, которые тесно связаны с распределением функции по граням. Кроме того, в случае $\lambda$-функции коэффициенты Фурье не равны нулю только в вершинах веса $h$. Поскольку значения $\lambda$-функции известны во всех вершинах веса не более $h$, то это, в частности, верно и для всех вершин $h$-мерной грани $\Gamma^{s(\alpha)}$ для произвольной вершины $\alpha\in W_h$. А значит, известно и $(s(\alpha),\beta)$-локальное распределение для произвольной вершины $\beta\in\Gamma^{s(\alpha)}$
(ясно, что для различных вершин с совпадающим носителем такие грани совпадут). По теореме о локальных энумераторах \ref{q-eig_loc} однозначно определены также $(\overline{s(\alpha)},\beta)$-локальные распределения, а потому и суммы значений
$$\eta(\alpha,\beta) = \sum_{\gamma\in\Gamma^{\overline{s(\alpha)}}(\beta)} f(\gamma)
 = \sum_{\gamma\in\Gamma^{\overline{s(\alpha)}}(\beta)} \varphi(\gamma)$$
в гранях $\Gamma^{\overline{s(\alpha)}}(\beta)$ для всех $\beta\in \Gamma^{s(\alpha)}$.
Вычислим коэффициент Фурье в произвольной вершине $\alpha\in W_h$:
$$\widehat{f}(\alpha)=
\sum_{\beta\in {\bf F}_q^n} f(\beta) \overline{\xi^{\langle\alpha,\gamma\rangle}}=
\sum_{\beta\in \Gamma^{s(\alpha)}} \sum _{\gamma\in\Gamma^{\overline{s(\alpha)}}(\beta)} f(\gamma) \overline{\xi^{\langle\alpha,\gamma\rangle}}=$$
$$\sum_{\beta\in \Gamma^{s(\alpha)}} \sum _{\gamma\in\Gamma^{\overline{s(\alpha)}}(\beta)} f(\gamma) \overline{\xi^{\langle\alpha,\beta\rangle}}=
\sum_{\beta\in \Gamma^{s(\alpha)}} \overline{\xi^{\langle\alpha,\beta\rangle}} \eta(\alpha,\beta) $$
Итак, все коэффициенты Фурье функции $f$ однозначно определяются по функции $\varphi$, а значит, и сама функция $f$ задается однозначно функцией $\varphi$.
\hfill $\square$

Соединяя эту лемму с теоремой \ref{main}, получаем:
\begin{theorem} \label{sphere-cub}
Пусть $\lambda$ -- собственное значение ${\bf F}_q^n$ с номером $h$  и $f$ -- $\lambda$-функция.
Пусть $\varphi: W_h\longrightarrow C$ -- произвольная функция, причем для всех $\alpha\in W_h$ выполнено $f(\alpha)=\varphi(\alpha) .$
Тогда 
все значения функции $f$ определяются однозначно функцией $\varphi$, если
$$\sum_{i=0}^k r^k_{i, d-k} P^{(q-1)}_i(l,k) \neq 0 , \ \ \ \ k=0,\ldots,h, \ \ l=0,\ldots,k,$$
где коэффициенты $r^k_{ij}$ определяются формулой (\ref{rij}) в случае $k\leq n-h$
и формулой (\ref{rijk>min}) в случае $k>n-h$.
\end{theorem}

Подмножество $L$ вершин гиперкуба называется \emph{ реконструктивным}, если из условий
$$f(\alpha)=0, \ \ \ \alpha\in L,$$
$$\widehat{f}(\alpha)=0, \ \ \ \alpha\in\overline{L}$$
следует, что $f$ тождественно равна нулю. Последняя теорема означает, что
при условиях (\ref{maincond}) сфера радиуса $d$ является реконструктивным множеством.

\end{document}